\documentclass[paper]{siamart171218}
\usepackage[english]{babel}
\usepackage{pgf}

\usepackage{tikz}
\usepackage[utf8]{inputenc}
\usepackage{amsmath,amssymb,amsfonts,todonotes}
\usepackage{enumerate}
\usepackage{blkarray}
\usepackage{cite}
\usepackage{mathtools}

\usetikzlibrary{matrix}

\usetikzlibrary{calc,fit,matrix,arrows,automata,positioning,patterns}

\usepackage{hyperref}
\usepackage{xcolor}
\hypersetup{
	colorlinks,
	linkcolor={blue!50!black},
	citecolor={blue!50!black},
	urlcolor={blue!80!black}
}
\newcommand{\B}[1]{\mbox{\boldmath $#1$}}

\usepackage{pgfplots}
\usepackage{hyperref}

\newenvironment{code1}{%
	\mathcode`\:="603A  % TeXbook pp 134, 154, 359 (top)
	\par
	\upshape
	\begin{list} % To give indentation
		{} {\leftmargin = 0.0cm}
		\item[]
		\begin{tabbing}
			\hspace*{.3in} \= \hspace*{.3in} \=
			\hspace*{.3in} \= \hspace*{.3in} \=
			\hspace*{.3in} \= \hspace*{.3in} \= \kill
		}{\end{tabbing}\end{list}}

\newcommand{\Rc}{\raisebox{-0.11cm}[0.1cm][0.0cm]{\mbox{$\Rsh$}}}
\newcommand{\rshno}{\rotatebox[origin=c]{180}{\reflectbox{\Rc}}}
\newcommand{\rc}{\raisebox{0.11cm}[0.1cm][0.0cm]{\mbox{\rshno}}}

\newcommand{\ca}{\hspace{-8.5pt}
	\makebox[0.3cm]{\raisebox{-0.2cm}[0.1cm][0.1cm]{{${\hookrightarrow}$}}}}

\newcommand{\Uk}{\mathcal{U}_k}

\newcommand{\stb}{\curvearrowright}

\newcommand{\slb}{\reflectbox{$\curvearrowright$}}

\newsavebox{\combarrow}
\savebox{\combarrow}[0.4cm]
{
	\Rc\hspace{-8.5pt}
	\raisebox{-0.4cm}{\rc}
}

\newcommand\TheTitle{Orthogonal iterations on structured pencils}
\newcommand\TheAuthors{R. Bevilacqua, G.M. Del Corso, L. Gemignani}

\title{\TheTitle\thanks{This work is partially supported by GNCS-INdAM.}}

\author{R. Bevilacqua\thanks{Dipartimento di Informatica, Universit\`a
		di Pisa, Pisa, Italy, \email{roberto.bevilacqua@unipi.it}} \and  G.M. Del Corso\thanks{Dipartimento di Informatica, Universit\`a
		di Pisa, Pisa, Italy, \email{gianna.delcorso@unipi.it}} \and   L. Gemignani\thanks{Dipartimento di Informatica, Universit\`a
		di Pisa, Pisa, Italy, \email{luca.gemignani@unipi.it}}}

\ifpdf
\hypersetup{
	pdftitle={\TheTitle},
	pdfauthor={\TheAuthors}
}
\fi

\begin{document}
	
	\maketitle
	
	\begin{abstract}
          We present a class of fast subspace tracking algorithms based on orthogonal iterations for structured matrices/pencils that can be represented
          as small rank perturbations of unitary matrices.   The algorithms rely upon  an updated data sparse factorization  --named LFR factorization-- using  orthogonal
          Hessenberg matrices. These  new subspace trackers reach a complexity of only $O(nk^2)$  operations per time update,
          where $n$ and $k$ are the size of the matrix and of the small rank perturbation, respectively.
        \end{abstract}
	
	\begin{keywords}
		Subspace tracking,  Orthogonal iteration, Eigenvalues, Eigenvectors, Unitary matrices, Low rank correction
	\end{keywords}
	
	\begin{AMS}
		65F15
	\end{AMS}
	\pagestyle{myheadings}
	\thispagestyle{plain}
	\markboth{R.~BEVILACQUA, G.~M.~DEL CORSO, L.~GEMIGNANI}{Orthogonal iterations on structured pencils}
	
        \section{Introduction}

Subspace tracking is an important  tool
in modern adaptive systems. The goal is the recursive
estimation of  the $s$
largest or smallest eigenvalues and the associated
eigenvectors of a possibly time-varying  matrix/pencil. In this paper we are concerned with
the design of fast subspace trackers for a certain class of structured matrices (pencils)  that can be represented as
small  rank perturbations of unitary matrices.  The  paramount example is (block) companion matrices/pencils.
  Applications arise in several different frameworks.
Subspace tracking is often required in signal
processing, especially in  multidimensional harmonic retrieval and system
identification algorithms because the characteristics of the signal can be  retrieved  from  the roots of some   associated
matrix polynomial \cite{Sinap,GuD}.  The computation of the roots of a matrix polynomial   also plays an important role  in the stability analysis of time-varying  dynamical systems,   which
amounts to establishing whether a linear system affected by some time-varying  parameters is asymptotically stable for all the admissible values of the parameters \cite{Galindo}.
More generally,   subspace tracking  is relevant for    model approximation and model reduction of dynamic systems \cite{Freund,Cherifi}.

Linearization and discretization are basic techniques  used for the search of a
reduced or approximated model.  As result, we usually determine a matrix polynomial which capture the dominant features of the original model.
A motivating  application   of this  approach    can be pursued for  solving  nonlinear eigenvalue problems  (NEP)  of the form $T(z)v=0,$ $v\neq 0$,
where $T:\Omega\to \mathbb{C}^{k\times k}$ is a holomorphic matrix-valued function and $\Omega \subseteq \mathbb{C}$ is a connected and open set.
The pair $(\lambda, v)$, $ v\neq 0$,  is an eigenpair of $T$ if it  satisfies $T(\lambda)v=0$, i.e., $\det(T(\lambda))=0$ and $v\in {\rm Ker}(T(\lambda))$.
Nonlinear eigenvalue problems arise in many applications~\cite{NLEVP,GT17}. The most studied case is the polynomial eigenvalue problem  (PEP) that can be tackled  by finding  suitable linearizations~\cite{MV04}
to  convert PEP into an equivalent generalized eigenproblem. Linearization methods using  block companion  forms ~\cite{HRT} allow the design of fast and stable methods~\cite{AMRVW17,BDG15,BDG2020}
which exploit the unitary plus low rank structure~\cite{AMRVW17}.  Various methods have been also  proposed in the literature to solve a NEP directly, for example Newton's method~\cite{Kr09}
and contour integrals~\cite{Be12,VK16,GMP18} techniques.  For a general   matrix-valued function $T(z)$  the associated eigenvalue problem
might have infinitely many eigenvalues,  and, hence, the usual scenario is to focus on computing  a few eigenvalues   located in a  certain subset $\Delta\subset \Omega$.
A possible approach   consists  first of  approximating  $T(z)$ with a matrix polynomial $P_\ell(z)$ inside $\Delta$,  and then of  computing the eigenvalues  of $P_\ell(z)$ to provide  numerical approximations of the eigenvalues of $T(z)$ in  $\Delta$. Since in general we are interested only in a few eigenvalues of $T(z)$, it is convenient to approximate only the
eigenvalues of interest rather than approximate the all spectrum of $P_\ell(z)$.

Subspace trackers based on orthogonal iterations  can be  numerically accurate and backward stable. The method of  orthogonal iteration  goes back to Bauer (see \cite{Rut} and the reference given therein).
If $s$ is the dimension of the subspace we want to approximate, a fast $O(ns)$ tracker for scalar companion matrices  based on orthogonal iteration  first appeared in \cite{STR}. The  algorithm is basically a game of  orthonormal Givens plane
rotations  moved from one side to the other side of orthogonal factors.   More recently these algorithms have been termed  core-chasing  ones \cite{R_book}.   In this paper  we extend the game to
more general  matrices $A\in \mathbb C^{n\times n}$  which are  unitary  plus some low rank-$k$ correction term. The development   follows by  exploiting the properties of a suitable LFR factorization  \cite{BDG2020,BDGHess20}
of  some  bordered
extension $\widehat A$ of  $A$, that is, $\widehat A=LFR$, where $L$ $(R)$ is a unitary $k$-lower ( $k$-upper) Hessenberg matrix  and $F=U+EZ$  is a unitary plus rank-$k$ matrix where   $U$  is a block diagonal
	unitary matrix of the form
	$\begin{bmatrix}I_k &\\ & \hat U \end{bmatrix} $ and  $E=[I_k, 0]^T $.  The unitary matrix $\hat U$ can be expressed as product of $\ell<n-k$ unitary Hessenberg matrices.
        It is shown that the shape of $U$ determines the shape of $\widehat A$. In particular $\widehat A$ and, a fortiori, $A$
        is upper triangular if and only if  $\hat U$ is upper triangular and hence diagonal.  Based on  this property it  makes possible to design an implementation of the inverse orthogonal iteration  scheme applied to
        $A$  using only  $O(n\max\{s,\ell\}k)$ ops per iteration.  The resulting algorithm is maximally fast w.r.t. the size of the correction term and it is  backward stable.   Moreover, it  can be easily generalized
        to deal  with both   the orthogonal iteration and the  inverse orthogonal iteration method  for structured pencils $(A,B)$ where $A$ and/or  $B$  are perturbed  unitary  matrices.

The paper is organized as follows.   In Section~\ref{two} we recall  the theoretical background concerning  the orthogonal iteration methods and
the properties of modified unitary matrices. Section~\ref{three} presents the derivation  of our  fast adaptations of the orthogonal iteration methods for modified unitary matrices.
In Section~\ref{four}
we show the results of numerical experiments that   lend   support   to   the   theoretical   findings. Finally, Section~\ref{five}  summarizes  conclusions and future work.

\section{Preliminaries}\label{two}
\hfill\\

In this section we recall  some preliminary results  concerning   the formulation of   both direct and inverse  orthogonal   iteration   schemes for matrix pencils  as well as   the  structural properties  and  data-sparse
representations
of modified unitary matrices.

\subsection{The Method of Orthogonal Iteration  for Matrix Pencils}\label{subsec1}
\hfill\\

The method of orthogonal iteration (sometimes  called subspace iteration
or simultaneous iteration) can be easily   generalized for matrix pencils.  Let $A-\lambda B$, $A,B\in \mathbb C^{n\times n}$, be a regular matrix  pencil with $A$ or $B$ invertible.  The orthogonal iteration
method can be  applied  for approximating the largest or smallest magnitude eigenvalues of   the matrix pencil by working on the  matrices $B^{-1} A$  or $A^{-1}B$.

If $B$ is nonsingular then  a generalization of the orthogonal  iteration method  to compute the $s$-largest (in magnitude) generalized eigenvalues  and corresponding eigenvectors of the matrix  pencil $A-\lambda B$
proceeds as follows:
\begin{equation}\label{orthit}
\left\{\begin{array}{ll}
AQ_i=BZ_{i+1} & i=1, 2 \ldots  \\
Q_{i+1}R_{i+1}=Z_{i+1} & \mbox{economy size QR factorization of } Z_{i+1}.
\end{array} \right.
\end{equation}
where $Q_1\in \mathbb C^{n\times s}$  is a starting orthonormal matrix comprising the initial approximations of the desired eigenvectors.  A detailed convergence analysis of this  iteration can be found  in \cite{Arb}.
It is found that  the convergence is properly understood in terms of invariant subspaces.  Specifically,  under mild assumptions it is proved that the angle between the subspace generated by the columns of $Q_i$ and
the invariant subspace associated with  the $s$ largest-magnitude generalized eigenvalues $\lambda_1, \lambda_2, \ldots, \lambda_s \in \mathbb{C}$, with
$|\lambda_1|\geq \ldots \geq |\lambda_s|>|\lambda_{s+1}| \geq \ldots \geq |\lambda_{n}|$, tends to zero as $O((|\lambda_{s+1}|/|\lambda_{s}|)^i)$   for $i$  going  to infinity.
An effective stopping criterion is  the following
\begin{equation} \label{eq:stop}
\|E_i\|=\|(I-Q_{i-1}Q_{i-1}^*)Q_i\|<\tau,
\end{equation}
where $\tau$ is the desired tolerance on the residual. Observe that this quantity measures the distance between the subspaces $\mathcal{S}_{i-1}=\rm{span}\{Q_{i-1}\}$ and $\mathcal{S}_{i}=\rm{span}\{Q_i\}$, in fact $(I-Q_{i-1}Q_{i-1}^*))Q_i$ can be taken as a measure of the angle between $\mathcal{S}_{i-1}$ and $\mathcal{S}_{i}$. Note moreover
that
$$
E_i^*E_i=Q_i^*(I-Q_{i-1}Q_{i-1}^*)(I-Q_{i-1}Q_{i-1}^*)Q_i=Q_i^*(I-Q_{i-1}Q_{i-1}^*)Q_i=I_s-W^*W,
$$
where $W=Q_{i-1}^*Q_i$. Then $\|E_i\|^2=1-\sigma_{min}^2(W)$. At convergence we expect $\mathcal{S}_{i-1}\equiv \mathcal{S}_{i}$ then $Q_{i}=Q_{i-1}U$
for an $s\times s$ unitary matrix $U$, then $\sigma_i(W)=\sigma_i(U)=1$ for $i=1, \ldots, s$.

Assume  now we are given a pencil $A-\lambda B$ with $A$ invertible  and that we would like to compute the $s$ smallest-magnitude generalized eigenvalues
$\lambda_1, \lambda_2, \ldots, \lambda_s \in \mathbb{C}$ with $|\lambda_1|\le |\lambda_2|\le \cdots\le |\lambda_s|<|\lambda_{s+1}| \leq \ldots \leq |\lambda_{n}|$.
Inverse orthogonal iterations can be used to approximate the desired eigenvalues. Starting with a set of $s$ orthogonal vectors
stored in matrix $Q_0\in \mathbb{C}^{n\times s}$, we compute the sequences
\begin{equation}\label{orthitinv}
\left\{\begin{array}{ll}
AZ_i=BQ_{i-1} & i=1, 2 \ldots  \\
Q_{i}R_{i}=Z_i & \mbox{economy size QR factorization of } Z_i.
\end{array} \right.
\end{equation}

Direct and inverse orthogonal iterations  \eqref{orthit},\eqref{orthitinv} can be carried out by solving the associated linear systems, but such an approach is prone to numerical instabilities due to the conditioning of
the resulting coefficient matrices.  A more accurate way to perform these schemes is using the QR factorization of the matrices involved. In particular  for the inverse iteration \eqref{orthitinv}
one may proceed at each step as follows:
\begin{enumerate}
\item Compute the full QR factorization of $Q_R R_R:=BQ_{i-1}$;
\item Compute the full RQ factorization of $R_LQ_L:=Q_R^*A$;
\item Determine  $Q_i$ such that $Z_i=Q_iR_i$ is the solution of  $R_LQ_L Z_i=R_R$. The set of orthogonal vectors satisfying the linear system is such that  $Q_i^*=Q_L(1:s, :)$.
\end{enumerate}

In the next subsection we introduce a suitable factorization of modified unitary matrices which makes possible to realize this  QR-based process in an efficient way.
Since  \eqref{orthit} can be implemented similarly by interchanging the role of the matrices $A$ and $B$ in the sequel we refer to orthogonal iteration as the scheme \eqref{orthitinv}.

\subsection{Fast compressed representations of modified unitary matrices}
\hfill\\

In this section we introduce a suitable compressed factorization  of unitary plus rank-$k$ matrices which can be exploited for
the design of fast  orthogonal iterations according to the QR-based process described above.
See~\cite{BDG2020,BDGHess20} for additional theoretical results on this factorization.

We denote by  $\mathcal{U}_k$ the set of unitary-plus-rank-$k$ matrices, that is, $A \in \mathcal{U}_k$  if and only if there exists a unitary matrix $V$
and two skinny matrices $X, Y\in \mathbb{C}^{n\times k}$ such that $A = V + XY^*$.
 A key role is played by generalized Hessenberg factors.
\begin{definition}\label{def:genhess}
	A matrix $R\in \mathbb C^{m\times m}$ is called {\em $k$-upper Hessenberg} if $r_{ij}=0$ when $i>j+k$.
	Similarly, $L$ is called  {\em $k$-lower Hessenberg} if $l_{ij}=0$ when $j>i+k$.   In addition, 	
	when  $R$ is {\em $k$-upper Hessenberg} ($L$ is {\em $k$-lower Hessenberg}) and
	the outermost entries are non-zero, that is, $r_{j+k,j}\neq 0$ ($l_{j,j+k}\neq 0$), $1\leq j\leq m-k$,
	then  the matrix is called {\em proper}.	A matrix which is simultaneously $k$-lower and $k$-upper Hessenberg is caller {\em $k$-banded}.
\end{definition}

Note that for $k=1$  a Hessenberg matrix  is proper if and only if it is unreduced.
Also,  a $k$-upper Hessenberg matrix $R\in \mathbb C^{m\times m}$ is proper if and only if $\det(R(k+1:m, 1:m-k))\neq 0$.
Similarly a $k$-lower Hessenberg matrix $L$ is proper if and only if  $\det(L(1:m-k, k+1:m))\neq 0$.
To make the presentation easier when possible, we use the letter $R$  to denote unitary generalized upper Hessenberg matrices, and the letter $L$ for unitary generalized lower Hessenberg matrices.

Note that $k$-lower (upper) Hessenberg matrices can be obtained as the product of $k$ matrices with the lower (upper) Hessenberg structure, and that
unitary block Hessenberg matrices with blocks of size $k$ are (non-proper) $k$-Hessenberg matrices.

In the following we will work with Givens rotations acting on two consecutive rows and columns. In particular we will denote by
${\mathcal G}_i=I_{i-1}\oplus G_i \oplus I_{n-i-1}$ the $n\times n$ unitary matrix where $G_i$ is a $2\times 2$ complex Givens rotation  of the form $\left[ \begin{array}{cc} c &-s\\s& \bar c \end{array}\right]$
such that $|c|^2+s^2=1$, with $s\in \mathbb{R}, s\ge 0$. The subscript index $i$ indicates the active part of the matrix
${\mathcal G}_i$. In the case  $G_i=I_2$ we say that $G_i$ is a trivial rotation.

\begin{definition}\label{def:datasparse}
Given a unitary matrix $U\in \mathbb{C}^{n\times n}$, we say that $U$ has a {\em data-sparse} representation if can be expressed as the product of $O({n})$ Givens matrices of the form ${\mathcal G}_i$ described before, possibly multiplied by a phase matrix.
\end{definition}
Note that the definition~\ref{def:datasparse} includes unitary (generalized) Hessenberg defined in~\ref{def:genhess}, CMV-matrices~\cite{CMV}, and other zig-zag pattern~\cite{Va11}, as well as the product of a constant number of these structures.

Next lemma shows how the product between data-sparse unitary terms can be factorized swapping the role of the two factors.

\begin{lemma}\label{lem:swap}
Let $R\in \mathbb{C}^{n\times n}$ be a unitary $k$-upper Hessenberg matrix and let $U$ a unitary matrix. Then there exist two unitary matrices $V$ and $S$ such that
$RU=VS$  where  $S$ is  $k$-upper Hessenberg and
$V=\begin{bmatrix}
I_k&\\ & \hat V.
\end{bmatrix}
$
Similarly, let $L$ be a unitary $k$-lower Hessenberg matrix and let $U$ be a unitary matrix.  Then there exist two unitary matrices $V$ and $M$ such that
$LU=VM$  where  $M$ is  $k$-lower Hessenberg and
$V=\begin{bmatrix}
\hat V &\\ & I_k.
\end{bmatrix}
$
\end{lemma}
\begin{proof}
	Let us partition $R$ and $U$ as follows
	$$
	R=\begin{bmatrix}
	R_{11}&R_{12}\\R_{21} & R_{22},
	\end{bmatrix}\quad
	U=\begin{bmatrix}
	U_{11}&U_{12}\\U_{21} & U_{22}
	\end{bmatrix}
	$$
	where $R_{21}$ is an upper triangular matrix of  size $n-k$, and $U_{11}$ is square of size $n-k$. Multiplying $R$ and $U$ and imposing the conditions on the blocks of the product $VS$ we get $S_{11}=R_{11}U_{11}+R_{12}U_{21}, S_{12}=R_{11}U_{22}$. Moreover, since $S$ should be a $k$-upper Hessenberg  matrix, we have that $S_{12}$ should be triangular. Hence $\hat V$ and $S_{21}$ can be computed as the Q and R factor of the QR factorization of $R_{21} U_{11}+R_{22}U_{21}$. Finally we set $S_{22}=\hat V^*(R_{21}U_{12}+R_{22}U_{22})$.
	Using the same technique we prove that there exist $V$ and $M$ such that $LU=VM$.
\end{proof}

\begin{definition}. The (lower) staircase of a matrix $A = (a_{i,j} ) \in \mathbb{C}^{n\times n}$ is the sequence $m_j(A),  1 \le j \le n$, defined as follows
$$
m_0(A)=0, \qquad m_j(A) =\max\{ m_{j-1}(A),      \max_{i>j}   \{i : a_{i,j} \neq 0 \}\}
$$
\end{definition}
The sequence $m_j(A)$ allows to represent the zero pattern of a matrix, in particular to identify zero sub-blocks in the matrix, in fact for each $1\le j\le n,$ it holds  $A(m_j(A)+1:n, 1:j)=0$
We note that proper $k$-upper Hessenberg matrices have $m_j(A)=j+k$ for $j=1, \ldots, n-k$, and $m_j(A)=n$, for $j=n-k+1, \ldots, n$.

\begin{lemma}\label{lemma:stair}
Let $A\in \mathbb{C}^{n\times n}$ be a matrix with staircase described by the sequence $\{m_j(A)\}$, for $1\le j\le n$ and let $T\in \mathbb{C}^{n\times n}$ be a non singular upper triangular matrix,  we have
$m_j(TA)=m_j(AT)=m_j(A)$ for $1\le j\le n$.
\end{lemma}
\begin{proof}
Let $B=TA$. We have $b_{ij}=\sum_{s=i}^n t_{is}a_{sj}$. Because of the staircase profile of $A$ we have
$a_{sj}=0$, for $s>m_j(A)$, hence $b_{ij}=0$ for $i>m_j(A)$, implying that $m_j(B)\le m_j(A)$.
To prove the equality of the staircase profile of $B$ and $A$ consider the entry $b_{mj(A), j}=t_{mj(A), m_j(A)} a_{m_j(A), j} $. If $a_{m_j(A), j}\neq 0$ we conclude that $m_j(B)=m_j(A)$, however it may happen that $a_{m_j(A), j}= 0$, but from the definition of staircase profile we knot that there exists an index $s, s<j$ such that $a_{m_j(A), s}\neq 0$, and $m_s(A)=m_j(A)$. Hence $b_{m_j(A), s}= r_{m_j(A), m_j(A)}a_{m_j(A), s}\neq 0$, implying that $m_j(B)=m_s(A)=m_j(A)$.
The proof that $m_j(AT)=m_j(A)$ can be carried on with a similar technique.
\end{proof}

Any unitary matrix of size $n$ can be factorized as the product of at most $n-1$ unitary upper Hessenberg matrices\footnote{The argument still holds if we take lower unitary Hessenberg in place of the upper Hessenberg matrices.}, that is $U=R_{n-1}R_{n-2}\ldots R_1\, D$, where each $R_{i}=\prod_{j=i}^{n-1}{\mathcal G}_j$ and $D$ is a phase matrix.
 To describe the representation and the algorithm we use a pictorial representation already  introduced  in several papers (compare
with~\cite{R_book} and the references given therein).  Specifically, the action of a Givens rotation acting on two consecutive rows of the matrix is depicted as
$
\begin{array}{c}
\Rc\\
\rc\end{array}
$.
Then a chain of ascending two-pointed arrows as below
$$
\begin{array}{c@{\hspace{1mm}}c@{\hspace{1mm}}c@{\hspace{1mm}}c@{\hspace{1mm}}c@{\hspace{1mm}}c@{\hspace{1mm}}c@{\hspace{1mm}}c@{\hspace{1mm}}c@{\hspace{1mm}}c@{\hspace{1mm}}c@{\hspace{1mm}}c@{\hspace{1mm}}c@{\hspace{1mm}}c@{\hspace{1mm}}c@{\hspace{1mm}}c@{\hspace{1mm}}c@{\hspace{1mm}}c}
	%{c|cccccccccccccc}
	                    &\Rc     &     &     &     &     &  &  &  &  \\
	                    & \rc    & \Rc    &     &     &  & &  &  &  \\
	                    &     &  \rc   &\Rc &  &  &     &  &  &  \\
	                    &     &     & \rc & \Rc &     &     &  &  &   \\
	                    &     &      &     & \rc    &\Rc     &     &  &  &   \\
	                    &     &      &     &     &  \rc   &\Rc     &  &  &    \\
	                    &      &     &     &     &     &   \rc  &\Rc  &  &    \\
	                    &     &     &     &     &     &           & \rc &  & \\	
\end{array}=  \begin{bmatrix}
\times & \times &\times &\times &\times & \times& \times & \times \\
\times & \times &\times &\times &\times & \times& \times & \times \\
&\times &\times &\times &\times & \times& \times & \times \\
&&\times &\times &\times & \times& \times & \times \\
&&&\times &\times & \times& \times & \times \\
&&&&\times &\times & \times& \times \\
&&&&&\times &\times & \times\\
&&&&&&\times & \times \\
\end{bmatrix}
= {\mathcal{G}}_1{\mathcal G}_{2}\cdots {\mathcal G}_{7}
$$
represents a unitary upper Hessenberg matrix (in the case of size 8 ). Some of the rotations may be identities (trivial rotations), and we might omit them in the picture. For example, in the above definition of $H_{i}$ we    only have non-trivial rotations  ${\mathcal G}_{j}$ for $j\ge i$, while the representation of $H_3$ is

$$
\begin{array}{c@{\hspace{1mm}}c@{\hspace{1mm}}c@{\hspace{1mm}}c@{\hspace{1mm}}c@{\hspace{1mm}}c@{\hspace{1mm}}c@{\hspace{1mm}}c@{\hspace{1mm}}c@{\hspace{1mm}}c@{\hspace{1mm}}c@{\hspace{1mm}}c@{\hspace{1mm}}c@{\hspace{1mm}}c@{\hspace{1mm}}c@{\hspace{1mm}}c@{\hspace{1mm}}c@{\hspace{1mm}}c}
%{c|cccccccccccccc}
&.     &     &     &     &     &  &  &  &  \\
&     & .    &     &     &  & &  &  &  \\
&     &     &\Rc &  &  &     &  &  &  \\
&     &     & \rc & \Rc &     &     &  &  &   \\
&     &      &     & \rc    &\Rc     &     &  &  &   \\
&     &      &     &     &  \rc   &\Rc     &  &  &    \\
&      &     &     &     &     &   \rc  &\Rc  &  &    \\
&     &     &     &     &     &           & \rc &  & \\	
\end{array}=  \begin{bmatrix}
1 &   &  &  &  &  &   &   \\
 & 1  &  &  &  &  &   &   \\
& &\times &\times &\times & \times& \times & \times \\
&&\times &\times &\times & \times& \times & \times \\
&&&\times &\times & \times& \times & \times \\
&&&&\times &\times & \times& \times \\
&&&&&\times &\times & \times\\
&&&&&&\times & \times \\
\end{bmatrix}
= {\mathcal{G}}_3{\mathcal G}_{4}\cdots {\mathcal G}_{7}
$$

Givens transformations can also interact with each other by means of the {\it fusion} or the {\it turnover} operations
(see~\cite{Raf_book}, pp.112-115).  The fusion operation will be depicted as
$
\begin{array}{cccc}
 \Rc\ca &\hspace{-0.3cm} \Rc  \\[-0.05cm]
 \rc & \hspace{-0.3cm}\rc  \\[-0.05cm]
\end{array}
%\mbox{  resulting in  }
%\begin{array}{cccc}
%& \Rc   \\[-0.05cm]
%&\rc  \\[-0.05cm]
%\end{array},
$
and consists of the concatenation of two Givens transformations acting on the same rows. The result is a Givens rotation multiplied possibly by a $2\times 2$ phase matrix.
The turnover operation allows to rearrange the order of some Givens trans\-for\-ma\-tions (see~\cite{Raf_book}).

Graphically we will depict this rearrangement of Givens transformations as follows:
\begin{equation*}
\begin{array}{ccccc}
\Rc & \stb&\Rc  \\[-0.05cm]
\rc & \Rc &\rc & \\[-0.05cm]
& \rc &     &  \\[-0.05cm]
\end{array}\quad
\to
\quad
\begin{array}{ccccc}
&     & \Rc &  \\[-0.05cm]
& \Rc & \rc & \Rc & \\[-0.05cm]
& \rc & & \rc &  \\[-0.05cm]
\end{array}\qquad
\mbox{  or }
\quad
\begin{array}{cccccc}
&   & \Rc &  \\[-0.05cm]
& \Rc & \rc & \Rc & \\[-0.05cm]
& \rc &$\rotatebox[origin=m]{180}{$\curvearrowleft$}$& \rc &  \\[-0.05cm]
\end{array}
\quad
\to
\quad
\begin{array}{ccccc}
\Rc &  &\Rc  \\[-0.05cm]
\rc & \Rc &\rc & \\[-0.05cm]
& \rc &     &  \\[-0.05cm]
\end{array}.
\end{equation*}
$$
\begin{array}{ccccc}
\Rc & \slb&\Rc  \\[-0.05cm]
\rc & \Rc &\rc & \\[-0.05cm]
& \rc &     &  \\[-0.05cm]
\end{array}\quad
\to
\quad
\begin{array}{ccccc}
&     & \Rc &  \\[-0.05cm]
& \Rc & \rc & \Rc & \\[-0.05cm]
& \rc & & \rc &  \\[-0.05cm]
\end{array}\qquad
\mbox{  or }
\quad
\begin{array}{cccccc}
&   & \Rc &  \\[-0.05cm]
& \Rc & \rc & \Rc & \\[-0.05cm]
& \rc &$\rotatebox[origin=m]{180}{$\curvearrowright$}$& \rc &  \\[-0.05cm]
\end{array}
\quad
\to
\quad
\begin{array}{ccccc}
\Rc &  &\Rc  \\[-0.05cm]
\rc & \Rc &\rc & \\[-0.05cm]
& \rc &     &  \\[-0.05cm]
\end{array}.
$$

Note that if the Givens transformations involved in turnover operations are all non-trivial also the resulting three new matrices are non trivial (see~\cite{AMVW15}).

In this paper we are interested in the computation of a few eigenvalues of matrices belonging to ${\mathcal{U}}_k$ by means of the orthogonal iteration schemes outlined in Subsection~\ref{subsec1}.
We will represent these matrices in the so-called {\em LFR format},  a factorization introduced in~\cite{BDG2020,BDGHess20}.

\begin{definition}\label{lfr}
We say that a matrix $A\in \mathbb C^{n\times n}, A\in \Uk$ is represented in the LFR format if  $(L,F,R)$ are matrices  such that:
	\begin{enumerate}
	\item $A=LFR$;
	%	\item $L\in \mathbb C^{n\times n}$ is the product of $k$ unitary lower Hessenberg matrices;
	\item $L\in \mathbb C^{n\times n}$ is a  unitary $k$-lower Hessenberg matrix;
	%	\item $R\in \mathbb C^{n\times n}$ is the product of $k$ unitary upper  Hessenberg matrices;
	\item $R\in \mathbb C^{n\times n}$ is a  unitary $k$-upper  Hessenberg matrix;
	\item $F=U+ E \,Z^*\in \mathbb C^{n\times n}$  is a unitary plus rank-$k$ matrix, where  $U$  is a block diagonal
	unitary matrix of the form
	$\begin{bmatrix}I_k &\\ & \hat U \end{bmatrix} $, with $\hat U$  unitary, $E=[I_k, 0]^T $
	and  $Z\in \mathbb C^{n\times k}$.
	\end{enumerate}
\end{definition}

Any matrix in $\Uk$ can be brought in the LFR format as follows. Let $A\in \Uk$, such that $A=V+XY^*$ , then  $L$ is  a $k$-lower unitary Hessenberg such that $L^*X=\begin{bmatrix}T_k\\0\end{bmatrix}$, where $T_k$ is upper triangular.
Then $A=L(L^*V+\begin{bmatrix}T_k\\0\end{bmatrix} Y^*)$. Using Lemma~\ref{lem:swap} we can rewrite $L^*V=UR$, where $R$ is unitary $k$-upper Hessenberg and $U=\left[\begin{array}{c|c}
I_k & \\ \hline & \hat U
  \end{array}\right]$ with $\hat U$ unitary. Bringing $R$ on the right we get our factorization, where $F=U+E Z^*$ and $Z=R\,Y\,T_k^*$.

The {\em LFR format} of modified unitary matrices is the key  tool  to developing fast and accurate adaptations  of the  orthogonal iteration schemes.  This  will be the subject of the next section.

\section{Fast Adaptations of the Orthogonal Iterations}\label{three}
\hfill\\

As underlined in Subsection~\ref{subsec1}, to compute the next orthogonal vectors
approximating a basis of the invariant subspace we need to compute the QR decomposition of $BQ_i$ and the RQ decomposition of $Q_R^*A$.

\subsection{RQ and QR factorization of unitary plus low rank matrices}\label{rqqr_sect}
\hfill\\

To recognize the triangular factors from the LFR decomposition of $A$ and $B$ it is useful to embed  the matrices of the pencil
into larger matrices obtained  edging the matrices with $k$ additional rows and columns.
Next theorem explains how we can get such a larger matrices still maintaining the unitary plus rank-$k$ structure.
\begin{theorem}\label{theo:rep}
	Let $A\in \mathbb{C}^{n\times n} , A\in \Uk$, then it is always  possible to construct a matrix of size $m=n+k$, $\hat A\in \Uk$ of size $m=n+k$ such that
	\begin{equation}\label{eq:bigA}
	\hat A=\begin{bmatrix}
	A& B\\ 0_{kn}& 0_{kk}
	\end{bmatrix} \quad \mbox{ for a suitable } B.
	\end{equation}
The unitary part of $\hat A$  can be described with additional $nk$ Givens rotations respect to the representation of the unitary part of $A$.
\end{theorem}
\begin{proof} Let $A=V+XY^*$, with $V$ unitary.  We assume that $Y\in \mathbb{C}^{n\times n}$ has orthogonal columns otherwise we compute the economy size QR factorization of $Y$ and we set $Y=Q$ and $X=XR^*$. Set $B=VY$, and consider the matrices
	\begin{equation} \label{eq:emb}
	\hat V=\begin{bmatrix}
	V-VYY^*& B\\Y^*&0_k
	\end{bmatrix}
	, \quad \hat X=\begin{bmatrix}
	X+B\\-I_k
	\end{bmatrix},\quad \hat Y= \begin{bmatrix}
	Y\\0_k
	\end{bmatrix}.
	\end{equation}
	We can prove that $\hat V$ is unitary by direct substitution. The last $k$ rows of  $\hat A=\hat V+\hat X\hat Y^*$ are null.

Matrix $V$ can be factorized as product of unitary factors which are related to the original players of $A$, namely $V, X$ and $Y$. In particular
$$
	\hat V=\begin{bmatrix}
	V&\\
	& I_k
\end{bmatrix}S\begin{bmatrix}-I_k&\\ &I_n\end{bmatrix} S^*,
$$
where $S$ is a $k$-lower Hessenberg matrix such that $$S^*\begin{bmatrix}Y\\-I_k\end{bmatrix}=\begin{bmatrix}\sqrt{2} I_k\\0\end{bmatrix}.$$
Such a $S$ always exists and is proper (see Lemma 3 in~\cite{BDG2020}).
\end{proof}

Note that the LFR format of $\hat A$ is such that $L$ is proper, since $\hat X$ has the last $k$ rows equal to $-I_k$ (see~\cite{BDG2020} Lemma 3).

\begin{theorem}\label{teo:rep1}
Let $L, R\in \mathbb{C}^{m\times m}$, $m=n+k$,
be two unitary matrices, where $L$ is a proper unitary
$k$-lower Hessenberg matrix and $R$ is a proper  unitary $k$-upper Hessenberg matrix.
Let $U$ be a block diagonal unitary  matrix of the form
$U=\left[\begin{array}{c|c}
I_k & \\ \hline & \hat U
\end{array}\right]$, with $\hat U$ $n\times n$ unitary. Let $F$ be the  unitary plus rank$-k$ matrix defined as $F=U+E\,Z^*$  with $ Z\in \mathbb C^{m\times k}$.
%  $T=[T_k, 0]^T$ and $T_k$  upper triangular and invertible.
Suppose that the matrix $\hat A=LFR$ satisfies the block structure in~\eqref{eq:bigA}.
%\begin{equation}\label{bs}
%\widehat A=\left[\begin{array}{cc} A & \ast \\ 0_{k, n} & 0_{k,k}\end{array}\right].
%\end{equation}
Then  $A=\hat A(1:n, 1:n)$ is nonsingular and has the same staircase profile as $\hat U$.\qed
\end{theorem}

\begin{proof}
Since $L$ is unitary, we have $L^*\hat A=FR$. Let partition $L$ and $R$ as follows $L=\begin{bmatrix} L_{11} &L_{12}\\L_{21}&L_{22}\end{bmatrix}$ with $L_{12}$ a $n\times n$  lower triangular matrix, and similarly partition $R$ in such a way $R_{21}$ is $n\times n$ upper triangular.
Then we get
$$
L_{12}^* A=\hat U R_{21}.
$$		
For Lemma~\ref{lemma:stair} we know that $L_{12}^* A$ has the same staircase profile as $A$, and $\hat U R_{21}$ has the same staircase profile as $\hat U$.\qed
\end{proof}
Next lemma helps us recognizing triangular matrices in the LFR format.
\begin{lemma} \label{lem:tri}
	If $\hat A=L(I+E \,Z^*)R$  satisfies the block structure in~\eqref{eq:bigA},  then $\hat A$ is upper triangular.
\end{lemma}		
\begin{proof}
We have $L_{21}^* A=R_{21}$. Because  $L$ is proper, the triangular block $L_{12}^*$ is nonsingular and  $A=({L_{12}^{*}})^{-1}R_{21}$. Hence $A$ is upper triangular because is the product of upper triangular factors. $\hat A$ is upper triangular as well because is obtained padding with zeros~\eqref{eq:bigA}.
\end{proof}

\begin{figure}
\begin{minipage}[c]{1\textwidth}
$$
\underbrace{\begin{array}{c@{\hspace{1mm}}c@{\hspace{1mm}}c@{\hspace{1mm}}c@{\hspace{1mm}}c@{\hspace{1mm}}c@{\hspace{1mm}}c@{\hspace{1mm}}c@{\hspace{1mm}}c@{\hspace{1mm}}c@{\hspace{1mm}}c@{\hspace{1mm}}c@{\hspace{1mm}}c@{\hspace{1mm}}c@{\hspace{1mm}}c@{\hspace{1mm}}c@{\hspace{1mm}}c@{\hspace{1mm}}c}
%{c|cccccccccccccc}
&&\Rc     &     &     &     &     &  &  &  &  \\
& \Rc& \rc    & \Rc    &     &     &  & &  &  &  \\
&     \rc   &    \Rc     &  \rc   &\Rc &  &  &     &  &  &  \\
&&    \rc   &    \Rc       & \rc & \Rc &     &     &  &  &   \\
&     &      &   \rc   &    \Rc     & \rc    &\Rc     &     &  &  &   \\
&     &      &     &    \rc   &    \Rc    &  \rc   &\Rc     &  &  &    \\
&      &     &     &     &    \rc   &    \Rc    &   \rc  &\Rc  &  &    \\
&     &     &     &     &              &   \rc  &            & \rc &  & \\	
\end{array}}_{R}
\underbrace{	\begin{array}{c@{\hspace{1mm}}c@{\hspace{1mm}}c@{\hspace{1mm}}c@{\hspace{1mm}}c@{\hspace{1mm}}c@{\hspace{1mm}}c@{\hspace{1mm}}c@{\hspace{1mm}}c@{\hspace{1mm}}c@{\hspace{1mm}}c@{\hspace{1mm}}c@{\hspace{1mm}}c@{\hspace{1mm}}c@{\hspace{1mm}}c@{\hspace{1mm}}c@{\hspace{1mm}}c@{\hspace{1mm}}c}
	%{c|cccccccccccccc}
	&& &  & &&\Rc     &     &     &     &     &  &    \\
&&	 &  & & \Rc& \rc    & \color{white}{\Rc}    &     &     &  & &  \\
&&	 &  & \Rc&     \rc   &    \Rc     &  \color{white}{\rc}   &\color{white}{\Rc} &  &  &     &   \\
&&	 &\Rc  & \rc&\Rc&    \rc   &    \color{white}{\Rc}       &  \color{white}{\rc} & \Rc &     &     &   \\
&&	 \color{white}{\Rc}& \rc &\Rc & \rc    &  \color{white}{\Rc}    &    \color{white}{\rc}  &   \color{white}{\Rc}     & \rc    &\Rc     &     &    \\
&\Rc&\color{white}{\rc} & \Rc & \rc&  \Rc   &  \color{white}{\rc}     & \color{white}{\Rc}    &    \color{white}{\rc}   &    \Rc    &  \rc   &\Rc     &   \\
\Rc&\rc&	\Rc & \rc &\color{white}{\Rc} &    {\rc}  &   \color{white}{\Rc}  &   \color{white}{\rc}  &   \color{white}{\Rc}  &    \rc   &    \color{white}{\Rc}    &   \rc  &\Rc      \\
\rc&&	\rc &  & \color{white}{\rc} &     &     \color{white}{\rc} &     &   \color{white}{\rc}   &              &    \color{white}{\rc}  &            & \rc  \\	
	\end{array}	}_{U}
%\hspace{-4.9cm}\begin{minipage}\resizebox{4.5cm}{3.5cm}{
%		\begin{tikzpicture}
%		%\draw[blue, very thin, dashed] (-2,-2)--(3, -3.5);
%		%\draw[red, very thin, solid] (-2, -2)--(-3.5, -3.5);
%		\draw[black] (-1,0)--(8, 0)--(3.5, 5)--cycle;
%		\end{tikzpicture}
%	}\end{minipage{}
	=	\underbrace{	\begin{array}{c@{\hspace{1mm}}c@{\hspace{1mm}}c@{\hspace{1mm}}c@{\hspace{1mm}}c@{\hspace{1mm}}c@{\hspace{1mm}}c@{\hspace{1mm}}c@{\hspace{1mm}}c@{\hspace{1mm}}c@{\hspace{1mm}}c@{\hspace{1mm}}c@{\hspace{1mm}}c@{\hspace{1mm}}c@{\hspace{1mm}}c@{\hspace{1mm}}c@{\hspace{1mm}}c@{\hspace{1mm}}c}
	%{c|cccccccccccccc}
	%&& &  & &&     &     &     &     &     &  &    \\
	&&	 &  & & &   &     &     &     &  & & \phantom{*}  \\
	&&	 &  & & &   &     &     &     &  & &  \\
	&& &  & &&\Rc     &     &     &     &     &  &    \\
	&&	 &  & & \Rc& \rc    & \color{white}{\Rc}    &     &     &  & &  \\
	&&	 &  & \Rc&     \rc   &    \Rc     &  \color{white}{\rc}   &\color{white}{\Rc} &  &  &     &   \\
	&&	 &\Rc  & \rc&\Rc&    \rc   &    \color{white}{\Rc}       &  \color{white}{\rc} & \Rc &     &     &   \\
	&&	 \color{white}{\Rc}& \rc &\Rc & \rc    &  \color{white}{\Rc}    &    \color{white}{\rc}  &   \color{white}{\Rc}     & \rc    &\Rc     &     &    \\
	&&\color{white}{\rc} &  & \rc&   &  \color{white}{\rc}     &   &    \color{white}{\rc}   &      &  \rc   &     &   \\
	%&\Rc&\color{gray}{\rc} & \Rc & \rc&  \Rc   &  \color{gray}{\rc}     & \color{gray}{\Rc}    &    \color{gray}{\rc}   &    \Rc    &  \rc   &\Rc     &   \\	
	%\Rc&\rc&	\Rc & \rc &\color{gray}{\Rc} &    {\rc}  &   \color{gray}{\Rc}  &   \color{gray}{\rc}  &   \color{gray}{\Rc}  &    \rc   &    \color{gray}{\Rc}    &   \rc  &\Rc      \\
	%\rc&&	\rc &  & \color{gray}{\rc} &     &     \color{gray}{\rc} &     &   \color{gray}{\rc}   &              &    \color{gray}{\rc}  &            & \rc  \\	
\end{array}\quad}_{V}
\underbrace{\hspace{-0.5cm}\begin{array}{c@{\hspace{1mm}}c@{\hspace{1mm}}c@{\hspace{1mm}}c@{\hspace{1mm}}c@{\hspace{1mm}}c@{\hspace{1mm}}c@{\hspace{1mm}}c@{\hspace{1mm}}c@{\hspace{1mm}}c@{\hspace{1mm}}c@{\hspace{1mm}}c@{\hspace{1mm}}c@{\hspace{1mm}}c@{\hspace{1mm}}c@{\hspace{1mm}}c@{\hspace{1mm}}c@{\hspace{1mm}}c}
%{c|cccccccccccccc}
&\Rc     &     &     &     &     &  &    \\
 \Rc& \rc    & \Rc    &     &     &  & &  &  \\
     \rc   &    \Rc     &  \rc   &\Rc &  &  &     &   \\
&    \rc   &    \Rc       & \rc & \Rc &     &     &    \\
     &      &   \rc   &    \Rc     & \rc    &\Rc     &     &    \\
     &      &     &    \rc   &    \Rc    &  \rc   &\Rc     &      \\
      &     &     &     &    \rc   &    \Rc    &   \rc  &\Rc     \\
     &     &     &     &              &   \rc  &            & \rc  \\	
\end{array}
}_{S}
	$$
	\end{minipage}
\hspace{-10.7cm}
\begin{minipage}[c]{0.5\textwidth}
\resizebox{4.65cm}{3.3cm}{
\begin{tikzpicture}
\draw[red] (0,0)--(4, 0)--(2, 3.3)-- cycle;
%\draw[blue] (0.7,1.10)--(3.3, 1.10);
\draw[blue] (0.7,1.10)--(3.3, 1.10)--(2, 3.3)-- cycle;

\end{tikzpicture}
}
\end{minipage}
\hspace{-2.3cm}
\begin{minipage}[c]{0.4\textwidth}
	\resizebox{4.65cm}{3.3cm}{
		\begin{tikzpicture}
		\draw[red] (0,0)--(4, 0)--(2, 3.3)-- cycle;
		\draw[blue] (0.12,0.2)--(2.82,0.2)--(1.47, 2.425)-- cycle;
		
		\end{tikzpicture}
	}
\end{minipage}
%\end{minipage}
\caption{An example of the swap Lemma~\ref{lem:swap}. Here $R$ is a unitary 2-upper Hessenberg matrix factorized as the product of two descending sequence of Givens rotations.  The unitary matrix $U$ is  represented in terms of a sparse set of rotations. When applying the rotations of $U$ to $R$ only the rotations in the blue triangle pop out (transformed by the turnover operations) on the left, while the remaining Givens transformations of $U$ are fused with the bottom transformations of $R$. In the picture we omit to represent a diagonal phase matrix which can be produced by the fusion operations.}\label{fig:1}
\end{figure}
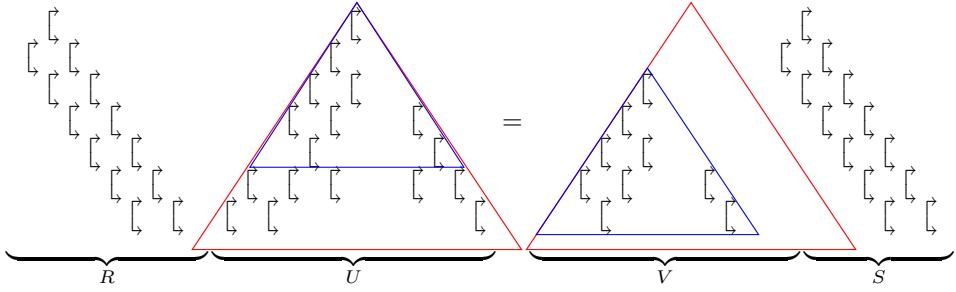

We now give an algorithmic interpretation of
Lemma~\ref{lem:swap}. A pictorial interpretation of  the lemma is given in Figure~\ref{fig:1}, where we omit the phase factors that are possibly present in the  general case.

Starting from Figure~\ref{fig:1} we can describe an algorithm for the ``swap'' of two unitary terms. In fact, we can obtain the new Givens rotations in the factors $V$ and $S$ simply applying repeatedly fusion and turnover operations as described by the  algorithm in Figure~\ref{alg:swap}. We formalize the algorithm as if the Givens involved in the swap were all non-trivial.
The algorithm has a cost $O(nk)$ only when $U$ admits a data-sparse representation. Matrix $U$ can be  generally  factorized as the product of at most $\ell\leq n-1$ unitary upper or lower Hessenberg matrices.
In this case the overall cost is $O(nk\ell)$. In procedure {\tt SwapRU} we choose to factorize $U$  as the product of lower Hessenberg factors, but we can obtain a similar algorithm expressing $U$ in terms of upper Hessenberg factors, and consider  the worst case $\ell=n-1$.
At step $i-1$ we have removed the first $i-1$ chains of ascending Givens rotations from $U$, and the situation is the following
$$
RU=V^{(1)}V^{(2)}\cdots V^{(i-1)} \, \Sigma_i L^{(i)}\cdots L^{(n-1)},
$$
where $\Sigma_i$ is an intermediate $k$-upper Hessenberg  which is transformed by the turnover and fusion operations. In particular $\Sigma_1=R$ and $\Sigma_{n-1}=S$.

At the step $i$ we pass the rotations in $L^{(i)}$, from right to left. The bottom $k$ Givens of each $L^{(i)}$ are fused with the Givens in the last rows of $\Sigma_i$, so that the shape of $\hat V$ reproduces the shape of the Givens rotations in the blue top triangle of $U$.

%\resizebox{8cm}{8cm}{
\begin{figure}\label{alg:swap}
\medskip
{\small{
	\framebox{\parbox{15.0cm}{
			\begin{code1}
				{\tt Procedure SwapRU }\\
				{\bf Input}:\quad   $R=R^{(k)}R^{(k-1)}\ldots R^{(1)}$, with $R^{(i)}={\mathcal G}^{(i)}_i {\mathcal G}^{(i)}_{i+1}\cdots {\mathcal G}^{(i)}_{n-1}$,\\
				\quad\quad $U=L^{(1)}\ldots L^{(\ell)}$, where $\ell\le n-1$ and $L^{(i)}=\Gamma^{(i)}_{n-1}\cdots \Gamma^{(i)}_{i}$\\
				\\
				{\bf for} $i=1, \ldots,k$, set $S^{(i)}=R^{(i)}$\\
			    {\bf for} $i=1:\ell$ \\
				\qquad {\bf for} $j=n-1:-1:n-k$ \{ {\em These Givens are fused} \}\\
				\qquad \quad {\bf if} $\Gamma^{(i)}_j\neq I_2$ \{ {\em only non trivial rotations are removed}\}\\
				\qquad \quad \quad after possible {\tt turnovers}, apply a {\tt fusion} with the last Givens in $S^{(n-j)} $\\
					\qquad \quad \quad
				\qquad \quad {\bf endif}\\
				\qquad {\bf endfor}\\
				\qquad {\bf for} $j=n-k-1:-1: i$ \\
				\qquad \quad {\bf if} $\Gamma^{(i)}_j\neq I_2$ \{ {\em only non trivial rotations are  removed}\}\\
				\qquad \quad \quad  $k$ turnover between Givens in $S$ acting on rows $j:n-1$ and $\Gamma^{(i)}_j$. \\
				\qquad \quad \quad The result is a Givens rotation $\tilde\Gamma_{j+k}^{(i)}$.\\
				\qquad \quad{\bf endif}\\
				\qquad {\bf endfor}\\
				{\bf endfor}\\
				{\bf Output}:\\
				\quad $S=S^{(k)}S^{(k-1)}\ldots S^{(1)}$, $V=\tilde L^{(1)}\ldots \tilde L^{(\ell)}$, where\\
				\quad $\tilde L^{(i)}=\tilde \Gamma^{(i)}_{n-1}\cdots \tilde \Gamma^{(i)}_{i+k}$.
			\end{code1}
		}
	}
}
%}
}
\caption{Procedure to swap  an upper and a lower generalized Hessenberg matrices.}
\end{figure}
\vspace{1cm}

Similarly to the procedure {\tt SwapRU}  we can design a procedure {\tt SwapLU} to factorize the product between a unitary $k$-lower Hessenberg matrix $L$ and a unitary matrix $U$ as the product of a unitary factor $V=\begin{bmatrix} \hat V & \\ & I_k \end{bmatrix}$ and a unitary $k$-lower Hessenberg matrix $M$. Note that from these two swapping procedures we can obtain also new factorizations when multiplying on the left a $k$-lower  or $k$-upper Hessenberg unitary matrix, that is  $U^*R^*=(RU)^*=(VS)^*=S^* V^*$, and $U^*L^*=(LU)^*=(VM)^*=M^*V^*$. We will denote the analog procedures as {\tt SwapUL} and {\tt SwapUR} keeping in mind that the matrices involved are unitary, and that we denote generalized lower Hessenberg matrices using the letter $L$ and generalized upper Hessenberg matrices using the letter $R$.

From the LFR format of $\hat A$ we can easily get the QR and RQ factorization of $\hat A$. This procedure requires $O(nk\ell)$ flops where $\ell$ is the number of Hessenberg unitary factors in $U$.
Let $\hat A=L(U+E \, Z^*)R$ as  in Definition~\ref{lfr}. Swapping  $L$ and $U$ according with Lemma~\ref{lem:swap}, i.e. $LU= Q \tilde L$, we have that $\tilde L(I+E\, Z^*)R$ is upper triangular for Lemma~\ref{lem:tri}.
Since $Q$ is unitary, we have a QR decomposition of  $A$.
 Similarly swapping $U$ and $R$ in such a way $UR=\hat  R \,Q$ we get an RQ decomposition of $A$ where the triangular factor is  $L(I+E\,\hat Z^*)\hat R$, with $\hat Z= UZ$, and the unitary factor is $Q$.
The proof is straightforward since, again from Lemma~\ref{lem:tri}, we have that $L(I+E\,,\hat Z^*)\hat R$ is upper triangular.

\subsection{The algorithm}
\hfill\\

In this section we describe the orthogonal iterations on a pencil $(A, B)$ where $A$ and $B$ are unitary-plus-low-rank matrices. For the sake of readability we  assume $A, B\in \Uk$, even if situation where the low-rank part of $A$ and $B$ do not have the same rank is possible: in that case we assume that $k$ is the maximum between the values of the rank part in $A$ and $B$.
%Typical examples ofe matrix pencils with this structure arise from the companion linearization of matrix polynomials, or from the  approximation of matrix functions on the unit circle and
%the subsequent Lagrange linearization~\cite{ACL}.

We  will assume that $A$ and $B$ have been embedded in larger pencil $(\hat A, \hat B)$ as described in Theorem~\ref{theo:rep}. Since $\det(\hat A-\lambda \hat B)=0$ for all $\lambda$, the pencil is singular and the $k$ new eigenvalues introduced with the embedding are indeterminate: Matlab returns ``NaN'' as eigenvalues in these cases.  However,  thanks to the block triangular structure of $\hat A$ and $\hat B$ the other eigenvalues coincide with those of the original pencil $(A, B)$. To guarantee that the orthogonal iterations on $(\hat A, \hat B)$ do not converge to an invariant subspace corresponding  with an indeterminate eigenvalue it is sufficient to start with a bunch of orthogonal vectors of the kind
$\hat Q_0=\begin{bmatrix} Q_0\\ 0_k \end{bmatrix}$. This guarantees that at each step we still have that $Q_i$ has the last $k$ rows equal
to zero and therefore the iterative process is basically applied to the smaller pencil.

We assume that $\hat A$ and $\hat B$ are in LFR format and we describe how we can carry out an  orthogonal iteration using only the turnover and fusion operations.
The key ingredient for the algorithm is the {\tt SwapUR} procedure and its variants as described  in Subsection~\ref{rqqr_sect}. In fact, working with the pencil and with the LFR factorization,
the orthogonal iterations can be reformulated as follows. Let
$$
\hat A=L_A(U_A+E,  Z_A^*) R_A, \quad \hat B=L_B(U_B+E\, Z_B^*) R_B
$$
be the LFR decomposition of $\hat A$ and $\hat B$.
The $s$ starting orthogonal vectors in $\hat Q_0\in \mathbb{C}^{N\times s}, N=n+k$, as well as all the intermediate orthogonal vectors $\hat Q_i$, can be represented as the product of $s$ sequences of ascending Givens rotations.  In fact, the columns of $\hat Q_0$  can be always be completed to an orthogonal basis $\{q_0, q_1, \ldots, q_N\}$ such that $[q_0, q_1, \cdots, q_N]$ is a $k$-lower Hessenberg matrix (see~\cite{BDG2020}).

\begin{figure} \label{algo:oi}
\centerline{
	\framebox{\parbox{8.0cm}{
			\begin{code1}
				{\tt Orthogonal Iterations }\\
				{\bf Input}:\quad LFR representation of $\hat A$ and $\hat B$, $\hat Q_0$ tolerance $\tau$, maxiter\\
				\\
				\ Represent $\hat Q_0$ with $s$ sequences of ascending Givens rotations\\	
				\ Using {\tt SwapUR}($U_A,R_A)$ compute the RQ decomposition of $\hat A=(L_A(I+E\, \hat Z_A^*)\hat R) Q_A$\\
				\ Using {\tt SwapLU}($L_B,U_B)$ compute the QR decomposition of $\hat B=Q_B(\tilde B(I+E\, Z_B^*)R_B)$\\		
				\ {\bf while} $i<$maxiter \& $\|E_i\|>\tau$ \\
				\qquad $Q_L:=${\tt MoveSequencesLeft}($B, \hat Q_i$)\\
				\qquad $Q_R:=${\tt MoveSequencesRight}($Q_L^*, A$)\\
				\qquad $\hat Q_{i+1}=Q_R$\\
				\qquad $E_i=(I-Q_{i-1}Q_{i-1}^*)Q_i$\\
				\qquad $i:=i+1$\\
				\ {\bf endwhile} \\
				\ $As=Q_i^*A Q_i$, $Bs=Q_i^*BQ_i$ \quad {\em $(As, Bs)$ is a pencil of dimension $s\times s$.}\\
				{\bf Output}:\\
				\ $eig(As, Bs)$
			\end{code1}
}
}}\caption{Inverse orthogonal iterations described ion terms of the LFR representation of the pencil.}
\end{figure}

\medskip
As we see in the algorithm in Figure~\ref{algo:oi}, the procedure boils down to the description of the two procedures {\tt MoveSequencesLeft} and {\tt MoveSequencesRight} that should be described in terms of the LFR representation.

	\medskip
\centerline{
	\framebox{\parbox{8.0cm}{
			\begin{code1}
				{\tt MoveSequencesLeft}\\
				{\bf Input}:\quad $\hat B=Q_B \tilde L_B(I+E\, \hat Z_B^*) R_B)$, $Q_0$ \\				
				\ $[ {\mathcal Q}_0,{ \mathcal R}_B]$={\tt SwapRU}$(R_B, Q_0)$\\
				\ $[{\mathcal P},  {\mathcal L}_B]$={\tt SwapLU}$(\tilde L_B,{\mathcal Q}_0 )	$\\	
				\ $[Q_L,  {\mathcal Q}_B]$={\tt SwapUL}$(Q_B,{\mathcal P})$\\				
				{\bf Output}: $Q_L$
			\end{code1}
		}
	}
}

	\medskip
\centerline{
	\framebox{\parbox{8.0cm}{
			\begin{code1}
				{\tt MoveSequencesRight}\\
				{\bf Input}:\quad $Q_L$, $\hat A=(L_A(I+E\, \hat Z_A^*)\hat R_A) Q_A$,\\
				\ $[ {\mathcal L}_A,{ \mathcal Q}_L^*]$={\tt SwapRU}$(Q_L^*, L_A)$\\
				\ $[{\mathcal R}_A, { \mathcal P}_L^*]$={\tt SwapUR}$({ \mathcal Q}_L^*,\hat R_A)	$\\	
				\ $[{\mathcal Q}_A, Q_R^*]$= {\tt SwapRU}$({ \mathcal P}_L^*, Q_A)	$\\			
				{\bf Output}: $Q_R$.
			\end{code1}
		}
	}
}

%\begin{gdc}
%Queste procedure sono un po' poco formali, mancano le diagonali di fase, e
%poichè $Q_L$ è $s$-lower Hessenberg anzichè applicare SwapUL in {\tt MoveSequenceRight}  sarebbe meglio applicare degli swap più precisi?
%\end{gdc}

\subsection{Measure of Backward stability} \label{Sec:bs}
\hfill\\

Suppose that our orthogonal iterations method  has reached a numerically  invariant subspace spanned by the $s$  orthogonal columns of matrix $Q$.
To measure the backward stability we analyze the quantity
$$
{\rm back_s}=\frac{\sqrt{2} \sigma_{s+1}([AQ, BQ])}{\|[A, B]\|}.
$$
This quantity is an upper bound to the usual backward stability measure.  Indeed,  we seek $\Delta_A$ and $\Delta_B$ such that
$$
(A+\Delta_A)Q=(B+\Delta_B)Q \Lambda,
$$
for a suitable invertible $\Lambda$. Since $A$ is invertible we may also suppose that $BQ$ is of maximum rank $s$.
We say that the algorithm is backward stable if
$$
\frac{\|[\Delta_A, \Delta_B]\|}{\|[A, B]\|}\approx O(\epsilon),
$$
where $\epsilon$ is the machine precision.

Let us consider the SVD decomposition of $[AQ, BQ]$. We expect that this matrix has $\sigma_{s+1}$ small since in  floating point
arithmetic $AQ=BQ \hat \Lambda$, with $\hat \Lambda$ invertible. We have
$$
[AQ, BQ]=[U1, U2]\begin{bmatrix}
\Sigma_1 &\\& \Sigma_2\\
\end{bmatrix}
\begin{bmatrix}
V_{11}^* & V_{21}^*\\ V_{12}^* &V_{22}^*
\end{bmatrix},
\quad \Sigma_1\in\mathbb{R}^{s\times s}.
$$
We find  $AQ=U_1 \Sigma_1 V_{11}^*+U_2\Sigma_2 V_{12}^*$ and $BQ=U_1 \Sigma_1 V_{21}^*+U_2\Sigma_2 V_{22}^*$.
Moreover the $s\times s$ matrix
$Q^*B^*BQ=V_{21}\Sigma_1^2 V_{21}^*+V_{22}\Sigma_2^2 V_{22}^*$ is invertible,
 then we have
$$
V_{21}\Sigma_1^2 V_{21}^*=Q^*B^*BQ(I-(Q^*(B^*B)Q)^{-1}V_{22}\Sigma_2^2V_{22}^*).
$$
Consider now the matrix
$I-(Q^*(B^*B)Q)^{-1}V_{22}\Sigma_2^2V_{22}^*$, which is invertible if
$$\|(Q^*(B^*B)Q)^{-1}V_{22}\Sigma_2^2V_{22}^*\|<1.
$$
This shows that  under this assumption $V_{21}$ is invertible as well.

Consider now the equality $\Delta_A Q-\Delta_B Q \Lambda=-A Q+BQ \Lambda$. Rewriting in terms of the SVD factors we get
$$
\Delta_A Q-\Delta_B Q \Lambda=U_1\Sigma_1(-V_{11}^*+V_{21}^*\Lambda)+U_2\Sigma_2(-V_{12}^*+V_{22}^*\Lambda).
$$
Since $V_{21}$ is invertible, among the infinite $s\times s$ matrices $\Lambda$  we can chose $\Lambda=V_{21}^{-*} V_{11}^*$, so that
$$
\Delta_A Q-\Delta_B Q \Lambda= -U_2\Sigma_2V_{12}^*+ U_2 \Sigma_2 V_{22}^*\Lambda.
$$
We can then set $\Delta_A=-U_2 \Sigma_2  V_{12}^*Q^*$ and  $\Delta_B=U_2 \Sigma_2 V_{22}^*Q^*$, and it holds
$\| \Delta_A\| \le \sigma_{s+1}([AQ, BQ])$ and $\| \Delta_B\|\le \sigma_{s+1}([AQ, BQ])$.
Finally, we conclude that
\begin{align*}
\|[ \Delta_A, \Delta_B]\|^2&=\rho(\Delta_A^* \Delta_A+\Delta_B^* \Delta_B)=
\\
&= \|\Delta_A^* \Delta_A+\Delta_B^* \Delta_B\|\le
(\|\Delta_A\|^2+\| \Delta_B\|^2)\le 2( \sigma_{s+1}([AQ, BQ]))^2.
\end{align*}

\section{Numerical results}\label{four}
\hfill\\

We perform several tests using nonlinear matrix functions $T(\lambda)\in \mathbb{C}^{k\times k}$. For matrix polynomials we consider the companion linearization while for non-polynomial matrix functions,
we first approximate the matrix function with polynomials of different degrees which are then linearized in pencils $(A, B)$ with $A, B\in {\mathcal U}_k$.

In all cases, when the pencil is built, our method performs the inverse orthogonal iterations  as defined in the Algorithm {\it Orthogonal Iterations} in Figure~\ref{algo:oi} until an invariant subspace is revealed.
Then the corresponding eigenvalues $\tilde \lambda_i$ are computed applying the Matlab
{\tt eig} function to the $s\times s$ pencil $(A_s,  B_s)$  determined  as the restriction of $A$ and $B$ to the subspace spanned by the columns of $Q_i$, i.e.
the generalized Rayleigh quotients of $A$ and $B$.

As a measure of the forward error we consider
$$
\mbox{err}_T(i)=\frac{\|T(\tilde \lambda_i) {\bf v}\|_2}{\|T(\tilde \lambda_i)\|_2 \|{\bf v}\|_2},
$$
where ${\bf v}$ is the $k$-th right singular vector of  $T(\tilde \lambda_i)$. In practice we compute
$$
\mbox{err}_T(i)= \frac{\sigma_k}{\sigma_1}
$$
where $\sigma_1\ge \sigma_2\ge \cdots \ge \sigma_k$ are the singular values of $T(\tilde \lambda_i)$.
 As a measure of convergence of the orthogonal iterations we consider instead
$$
\mbox {averr}_p=\frac{1}{s}\sum_{i=1}^s |\tilde \lambda_i-\mu_i|,
$$
where $\mu_i$ are the  ``exact'' eigenvalues of the pencil $(A, B)$ obtained with Matlab {\tt eig}.
Note that in the general non-linear case $\mbox {averr}_p$ refers to the average error with respect to the zeros of the  approximating polynomial while $\mbox{err}_T(i)$ is the total error with respect to the zeros of the  nonlinear function $T(\lambda)$.
Hence,  small values of $\mbox {averr}_p$ guarantee the good behavior of the inverse orthogonal iterations on the pencil, while $\mbox{err}_T(i)$ measures also the quality of the approximation of the nonlinear
function with the matrix polynomial.

We tested our method  on some  matrix polynomials of  degree greater than 2 from  the NLEVP collection~\cite{NLEVP}  using the companion linearization.
Table~\ref{Tab:NLEP} summarizes the results. For the polynomial {\tt plasma\_drift} we repeated the experiment twice by setting the  stopping criteria in~\eqref{eq:stop} to a tolerance of $1.0e-4$,
and the maximum number of iterations to $450$, respectively.  This is  a very challenging problem for any eigensolver \cite{HP20}  with several eigenvalues  of high  multiplicity  and/or
clustered around zero.
In this case the estimated backward error is not so significant since our measure of backward stability assumes that we are at convergence and we have identified an invariant subspace.
In~\cite{HP20} the authors proposed a variation of the Jacobi-Davidson method for computing several eigenpairs of  the polynomial eigenvalue problem.
For the {\tt plasma\_drift} problem they fixed a  residual threshold of $1.0e-2$ and within  200 iterations
they were able to compute the approximations  of the 19 eigenvalues closer to the origin. The   performance of our method for  the  approximation of the same set of eigenvalues is examined  in Table~\ref{Tab:NLEP}.
%As reported in that paper the forward error is only slightly lower then $1.0e-2$.
Regarding  the other tests, by comparison of our  results  with those reported in ~\cite{AMRVW17} where a structured version of the QZ method is employed, we see that  for the {\tt orr\_sommerfeld} problem  we get a higher backward stability. However,  we are estimating only 2 or 4 eigenvalues while the QZ allows  to approximate all the spectrum
and, moreover,  differently from ~\cite{AMRVW17}  our error analysis assumes  an uniform  bound for the norm of the perturbation of $A$ and $B$.  The accuracy of the computed eigenvalues
is in accordance with the conditioning estimates. 
Our method on the {\tt orr\_sommerfeld} problem has a backward stability similar to that reported in~\cite{SSRR} where a balanced version of the Sakurai-Sugiura method with Rayleigh-Ritz projection was presented. For the {\tt butterfly} problem  our method achieves a higher backward stability w.r.t. ~\cite{SSRR}.
The number of iterations in {\tt relative\_pose\_5pt}  agrees with the separation ratio of the eigenvalues.
For the other tests there are remarkable differences in the number of iterations depending on the sensibility of
our stopping criterion \eqref{eq:stop} used in  Algorithm {\it Orthogonal Iterations} w.r.t.  specific features of the considered eigenproblem.  Comparisons with other
stopping criteria introduced in the literature is an ongoing work.

\begin{table}
	\small{
		\begin{center}
			\begin{tabular}{|l|l|l|l|l|l|l|l|l|}
				\hline
				name  & deg &k & s& ${|\lambda_s|}/{|\lambda_{s+1}|}$  & ${\rm averr_p}$ & it & ${\rm back_s}$ \\ \hline
				%				CL  & 16  & 5.19e-01         & 4.40e-03                          & butterfly & c & 4  & 64&  4& 1.70e-01 & 7.76e+05 &  1.71e-14 &1.25e-14 & 48 &3.29e-15 &2.15e-15 &4.80e+04 &3.00e+03\
				%butterfly &  4& 64& 8& 9.84e-01 & 7.74e+01 &  3.86e-08 &1.16e-07 & 984 &3.40e-08 \\ \hline
				{\tt butterfly} &  4& 64& 4& 9.53e-01 &5.45e-14 & 654 &5.01e-15  \\ \hline
				{\tt orr\_sommerfeld} &  4& 64& 2& 9.95e-01  &5.75e-06 & 22 &1.92e-18 \\ \hline
				{\tt orr\_sommerfeld} &  4& 64& 4& 9.91e-01  &  8.33e-06 & 28 &1.92e-18  \\ \hline
				{\tt plasma\_drift} & 3& 128& 19& 9.98e-01  &6.60e-02 & 6 &7.38e-06 \\ \hline
				{\tt plasma\_drift} & 3& 128& 19& 9.98e-01 &5.85e-04 & 450 &7.92e-08 \\ \hline	
				{\tt relative\_pose\_5pt} &  3& 10& 4& 3.20e-01  &1.99e-14 & 33 &1.34e-15\\ \hline
				
			%	planar\_waveguide & 4& 129& 4& 9.06e-01 & 1.48e+03  &3.52e-13 & 281 &5.21e-15 \\ \hline
				
			\end{tabular}
		\end{center}
	} \caption{ Matrix polynomials of degree 3 and 4 from the NLEVP collection.\label{Tab:NLEP}}
\end{table}

Another set of experiments   have dealt  with root-finding  for  a nonlinear matrix function.
Consider the holomorphic nonlinear matrix-valued function $T:\Omega\to \mathbb{C}^{k\times k}$, with $\Omega\subset \mathbb{C}$ a connected and open set and let $\lambda_i$ be the
$i-th$ exact eigenvalue of $T$, that is $T(\lambda_i)\B {v}=0$, with  $\B v$ a corresponding eigenvector. Computing an approximation of $\lambda_i$ can be accomplished as follows.
We first approximate the nonlinear function with a matrix polynomial of a given degree $d$ interpolating on the roots of unity or on the  Chebyshev points. With these choices
of points we have theoretical results~\cite{BDG18b,EK} about the uniform convergence of the interpolating polynomials to the nonlinear function inside the unit disk.
Let $P_d(z)$ be the approximating polynomial of degree $d$, then we  may consider  suitable
linearizations which give us a unitary plus-low-rank pencil $(A, B)$. The eigenvalues of this pencil provide an approximation of the  zeros of $T(z)$ inside  $\Omega$.

Specifically,  we tested  our algorithm on the companion linearization,
and on the unit diagonal plus-low-rank linearization obtained from the companion linearization applying a block-Fourier transform which diagonalizes the unitary part of $A$ (see~\cite{BDG18b} for more details).
The same structure  can also  be obtained starting from an ``arrowed linearization'' similar to the one proposed in~\cite{ACL} where the interpolating polynomial is written in the Lagrange basis.
Combining the different linearizations with the different choices of the nodes (roots of unity, roots of unity plus the origin, Chebyshev nodes) the following cases are treated:
%Vogliamo provare la arrowed con le radici ruotate come in 1.52?
%Vogliamo provare anche la arrower 2k? (qui occorre mettere mano al programma)
%
%In particular we considered the following linearizations
\begin{itemize}
	\item Companion linearization on the roots of unity (denoted in the tables by ``CL'')
	\item Companion linearization on the roots of unity plus the origin (denoted in the tables by ``CL0')
	\item Companion linearization on the Chebyshev roots (denoted in the tables by ``T')
	\item Diagonal linearization on the roots of unity (denoted in the tables by ``DL'')
\end{itemize}
We note that in the diagonal linearization the diagonal factor contains the interpolation nodes, hence to guarantee that this factor is unitary we can  choose as interpolation nodes only the roots of unity but not the other choices of nodes considered in the companion linearization such as the Chebyshev points or the origin.

We tested several non-linear matrix-valued functions found in the literature.
%NLEVP collection~\cite{NLEVP} as well as a function $F(z)$ from~\cite{AST}. In \cite{AST}the authors
%used techniques based on contour integrals to compute the eigenvalues of non linear functions. The functions tested were eventually shifted and scaled so that the eigenvalues will lie inside the unit disk.
Below is a description of these  functions.

\noindent{\bf Time-delay equation} \cite{Be12}.  The matrix function is  $T(z)=z+T_0+T_1 \exp(6z-1)$
with
$$
T_0 =\left[\begin{array}{cc}4 & -1\\-2 &5\end{array}\right]; \quad
T_1 = \left[\begin{array}{cc} -2 & 1\\ 4 &-1\\ \end{array}\right].
$$
This function has three eigenvalues inside the unit circle.

\noindent{\bf Model of cancer growth} \cite{BKZ12}. The matrix function is  $T(z)= z-A_0-A_1\exp(-r z)$, where
%	k=3; m=32;
%	r=5; b1=0.13; bq=0.2; mu1=0.28; mu0=0.11; muq=0.02; mug=0.0001;
%
%{\scriptsize{
		$$	A_0=\left[\begin{array}{ccc}-\mu_1& 0&0\\ 2b1& -\mu_2& b_Q\\ 0& \mu_Q& -(b_Q+\mu_G)\end{array}\right],\,
		A_1=\exp(-\mu_2r) \left[\begin{array}{ccc}2b_1& 0& b_Q\\ -2b_1& 0& -b_Q\\ 0&0&0\end{array}\right].
		$$
%}}
The parameters are chosen as suggested in~\cite{BKZ12}   by
setting $r=5; b_1=0.13; b_Q=0.2; \mu_1=0.28; \mu_0=0.11; \mu_Q=0.02; \mu_G=0.0001$, $\mu_2=\mu_0+\mu_Q$.
We refer  to~\cite{BKZ12} for the physical meaning of the constants and for the description of the model. This function has three eigenvalues inside the unit circle.

\noindent{\bf Neutral functional differential equation} \cite{ERL}.
The function is  scalar  $t(z)= -1+0.5 z+z^2+ h z^2 \exp(\tau z)$.
The case $ h=-0.82465048736655, \tau=6.74469732735569$
is analyzed in \cite{Kra} corresponding to a
Hopf bifurcation point.  This function has three eigenvalues inside the unit circle.

\noindent{\bf  Spectral abscissa  optimization} \cite{MBN}.  The function is
$T(z)=zI_3  -A -B \exp(-z \tau)$ with $\tau=5$, $B=b q^T$ and
$$
A=\left[\begin{array}{ccc}-0.08 &  -0.03&  0.2\\
0.2 &  -0.04&  -0.005\\  -0.06 &  0.2 &  -0.07
\end{array}\right], \ b=\left[\begin{array}{c}
-0.1\\ -0.2\\ 0.1
\end{array}\right], \
q=\left[\begin{array}{c} 0.47121273 \\ 0.50372106 \\0.60231834
\end{array}\right].
$$
Abscissa  optimization   techniques
favor multiple roots  and  clustered eigenvalues with potential numerical difficulties. This function has 4 eigenvalues inside the unit circle.
%({|bf capire perche abbiamo preso s=2})

\noindent{\bf Hadeler  problem} \cite{NLEVP}.  The matrix function is $T(z)=(\exp(z)-1)A_2+z^2 A_1-\alpha A_0$ where
$A_0, A_1,A_2\in \mathbb{R}^{k\times k}$, and $V=\mbox{ones}(k, 1)*[1:k]$
$A_0=\alpha I_k, A_1= k*I_k + 1./(V+V'), A_2=(k+1 - \max(V,V')) .* ([1:k]*[1:k]')$.
In our experiments we set $k=8$ and $\alpha=100$.
Ruhe~\cite{Ru73} proved that the problem has $k$ real and positive eigenvalues,
in particular two of them are $0<\lambda<1$, and hence lie inside the unit circle.

\noindent{\bf Vibrating string} \cite{NLEVP,So06}.
The model refers to a string of unit length clamped at one end, while the other one is
free but  is loaded with a mass $m$ attached by an elastic spring of stiffness $k_p$. Assuming $m=1$,
and discretizing the  differential equation one gets
the non linear eigenvalue problem $F(z)v=0$, where  $F(z)=A-B z+k_p C \frac{z}{z-k_p}$ is rational,
$A, B, C\in \mathbb{R}^{k\times k}$,
$k_p=0.01, h=1/k$,
{\scriptsize{
		$$
		A=\frac{1}{h} \left[\begin{array}{ccccc}
		2 & -1 &&&\\
		-1 &\ddots& \ddots &&\\
		& \ddots& & 2 &-1\\
		& & & -1 &1
		\end{array}\right],   \,
		B=\frac{h}{6} \left[\begin{array}{ccccc}
		4 & 1 &&&\\
		1 &\ddots& \ddots &&\\
		& \ddots& & 4 &1\\
		& & & 1 &2
		\end{array}\right], \,
		C=e_ke_k^T.
		$$
}}

\noindent{\bf The function $F(z)$} $: \Omega\to \mathbb{C}^{3\times 3}$ {\bf from}~\cite{AST} is defined as follows:
\begin{equation} \label{Fz}
F(z)=\left[\begin{matrix}
 2 e^z+cos(z)-14& (z^2-1)\sin(z)+(2 e^z+14)\cos(z)& 2 e^z-14 \\
(z+3)(e^z-7)&\sin(z)+(z+3_(e^z-7)\cos(z)& (z+3)(e^z-7) \\
 e^z-7&(e^z-7)\cos(z)& e^z-7
\end{matrix}\right].
\end{equation}
This function has six real known eigenvalues  $\left\{\pm \pi, \pm \pi/2, 0, \log(7) \right\}$.
We applied the transformation $z\to 4z+1$ to bring five of the six eigenvalues inside the unit disk. With this transformation we do not get an approximation of the eigenvalue $-\pi$ which after the translation is not inside the unit disk.

For all these problems we computed the interpolating polynomials over the roots of unity, the roots of unity plus the origin or the Chebyshev points in the range $[-1, 1]$ of different degrees, and we compared the performance of our algorithm on the different linearizations.  Linearizations based on interpolating at the roots of unity (plus the origin)  perform very similarly with
negligible differences in convergence and accuracy.  The linearization using the Chebyshev points can suffer of numerical instabilities for large degrees of the interpolating polynomial due to  the transformation of the polynomial basis.   For the sake of brevity  we report  here only the best result obtained using the polynomial with lower degree which guarantees the best performance in terms of forward error $\mbox{err}_T(i)$.
Ask  the authors for the complete set of results. In Table~\ref{Tab:sum} are summarized the results for the nonlinear functions considered.
In Table~\ref{Tab:Fz} are reported the complete  results for the Function in~\cite{AST}.

\begin{table}
	\small{
		\begin{center}
			\begin{tabular}{|l|l|l|l|l|l|l|l|l|}
				\hline
				 name&m(deg) &s&${|\lambda_s|}/{|\lambda_{s+1}|}$ & ${\rm err_T(1)}$ &  ${\rm err_T(s)}$  & ${\rm averr_p}$ & it & ${\rm back_s}$ \\ \hline
				Time-del& CL(64)  & 3& 5.25e-01                      & 5.74e-16         & 2.91e-14         & 7.56e-14        & 45 & 1.95e-15       \\ \hline
				Cancer&CL(32)  &  2& 5.18e-01                     & 4.22e-16 & 4.78e-15 & 5.76e-15 & 48 & 1.91e-15       \\ \hline
			    Neutral& CL0(64) & 2 &4.02e-01 						  & 6.72e-13 & 9.86e-13 & 1.63e-13 & 33  & 5.33e-16 \\ \hline
			    Spec.-abs& T(32) &  4& 1.02e-01  &                  1.66e-16 &1.16e-16   &2.12e-09 &24 &4.34e-16 \\ \hline
			    Hadeler &DL(32) &2& 6.35e-01  & 5.58e-16& 1.65e-14&7.30e-15&87 &3.14e-15 \\  \hline
			    Vib-str. & CL0(32) & 1& 5.10e-01& 5.24e-15& -&4.27e-14&64 &4.64e-16 \\ \hline
					\end{tabular}
		\end{center}
	} \caption{Best results for the non linear non-polynomial matrix functions.\label{Tab:sum}}
\end{table}

%In Tables~\ref{Tab:TD},\ref{Tab:CG},\ref{Tab:NF},\ref{Tab:SA},\ref{Tab:H},\ref{Tab:VS} and \ref{Tab:Fz} are reported the results obtained in our experiments.

We underline that the algorithm, accordingly with the measure of backward stability in Section~\ref{Sec:bs}, behaves as a backward stable method in every case. The number of iterations needed to meet the stopping condition which was set to $1.0e-14$ reflects the rate of convergence of the orthogonal iterations which depends on the ratio $|\lambda_s|/|\lambda_{s+1}|$. The values $\mbox{averr}_p$ measure the effectiveness of orthogonal iterations to approximate the eigenvalues of the pencil and of course are affected by the conditioning of the problem.
When Chebyshev points are used as interpolation nodes sometimes the pencil obtained is seriously ill conditioned and hence both Matlab {\tt eig} and our algorithm return inaccurate results.
%Note that also the process of retrieving the coefficients of the interpolating polynomial on the Chebyshev points is an operation prone to instability because, for large values of $d$ we need to solve the Vandermonde system which can be ill conditioned. This reflects on both the error $\mbox{err}_T$ and $\mbox{averr}_p$. For example for the  {\bf Neutral function} we could not retrieve the coefficients in the monomial basis of the Chebyshev interpolant of degree 128.
On the contrary, when working with the roots of unity the coefficients of the interpolating polynomial in the monomial basis are computed by means of an FFT which is very stable.

In general, for sufficiently large  values of the degree we get a very good approximation of the eigenvalues inside the unit disk.
We can compare these results with those reported in~\cite{BDG18b} and~\cite{AMRVW17} where respectively a QR and QZ method were employed to compute all the eigenvalues of the matrix/pencil. We see that
our results are comparable with those obtained in the literature but we need less operations. In fact the algorithms based on QR or QZ need $O(\mbox{d}^2 \,k^3)$ flops, while here we need $O(\mbox{d} \, k^2 \,s\, \mbox{it})$ flops. In general the number of iterations does not depend on the size of the problem, but only by on the ratio $|\lambda_s|/|\lambda_{s+1}|$, so the cost of the orthogonal iterations  can be asymptotically
lower and we do not have any advantage in computing all the eigenvalues since only $s$ of them are reliable because the polynomial is a good approximation of the non-linear function only inside the unit disk.
Comparing the results in Table~\ref{Tab:Fz} with those reported in paper~\cite{AST} we see that our results are much better, in particular  when using the method CL, that is the approximation of the non linear function with the interpolating polynomial over the roots of unity  combined with the companion linearization. In particular using the same degree ($\mbox{d}=64$) as in \cite{AST} we get a results with 5 more digits of precision respect to the results reported in~\cite{AST}.

\noindent
\begin{table}
	\small{
		\begin{center}
			\begin{tabular}{|l|l|l|l|l|l|l|l|l|l|l|} \hline
							m & deg &$\displaystyle\frac{|\lambda_s|}{|\lambda_{s+1}|}$& ${\rm err_T(1)}$   & ${\rm err_T(2)}$  & ${\rm err_T(3)}$  & ${\rm averr_p}$   & it & ${\rm back_5}$ \\ \hline
				CL  & 16 & 6.77e-01 & 3.68e-07 & 1.52e-06 & 3.57e-06 & 1.32e-12 & 81  & 3.23e-15 \\ 		\hline
				CL0 & 16 & 6.77e-01  & 5.25e-08 & 1.89e-07 & 2.66e-06 & 4.14e-12 & 82  & 4.31e-15 \\		\hline
				T   & 16 &7.10e-01 & 2.99e-12 & 2.29e-10 & 2.64e-10 & 1.21e-12 & 86  & 3.89e-15 \\		\hline
				DL  & 16 &6.92e-01 & 5.25e-08 & 3.62e-07 & 2.67e-06 & 1.10e-13 & 131 & 1.79e-15 \\		\hline
				CL  & 32 &6.92e-01 & 1.76e-16 & 5.40e-16 & 5.24e-15 & 4.62e-12 & 86  & 3.87e-15 \\		\hline
				CL0 & 32 &6.92e-01 & 1.69e-16 & 2.61e-16 & 1.63e-15 & 4.26e-12 & 82  & 4.18e-15 \\ 		\hline
				T   & 32 &6.92e-01 & 6.87e-16 & 2.85e-15 & 5.85e-14 & 1.94e-12 & 83  & 4.06e-15 \\		\hline
				DL  & 32 &6.92e-01 & 9.74e-16 & 1.69e-15 & 2.19e-15 & 1.08e-13 & 101 & 2.95e-15 \\		\hline
				CL  & 64 &6.92e-01 & {\bf 2.59e-17} & {\bf 4.26e-17 }& {\bf 1.14e-16 }& 2.65e-11 & 83  & 3.31e-15 \\ 		\hline
				CL0 & 64 &6.92e-01 & 2.71e-17 & 3.03e-17 & 1.95e-16 & 7.61e-12 & 82  & 3.69e-15 \\		\hline
				DL  & 64 &6.92e-01 & 5.40e-15 & 6.25e-15 & 8.61e-15 & 1.56e-13 & 101 & 2.45e-15 \\ \hline
			\end{tabular}
\end{center}
        } \caption{{\bf Function $F(z)$ in~\eqref{Fz}}.  For the inverse orthogonal iterations we used $s=5$, but we show the results only for the first three eigenvalues (corresponding to the value $\pi/2, \log(7)$ and 0).
         The other  remaining two roots are  approximated just as  well.   \label{Tab:Fz}}
\end{table}

\section{Conclusions and Future Work}\label{five}
\hfill\\

In this paper we have presented a fast and backward stable  subspace tracker for block companion forms using orthogonal iterations.
The proposed method exploits the properties of a  suitable data-sparse factorization of the matrix involving unitary factors.  The method can be
extended to more generally perturbed unitary matrices  and it can incorporate the acceleration techniques  based on the  updated computation of Ritz eigenvalues and eigenvectors
\cite{Arb}.  The design of  fast adaptations using  adaptive shifting techniques  such as the ones proposed in \cite{Jung} is an ongoing research project.

\bibliographystyle{siamplain}
\bibliography{compimpl}

	\end{document}